\newtheorem{thm}{Theorem}

\newtheorem{pred}{Proposition}

\newtheorem{gip}{Conjecture}

\newcommand{\dvo}{{\it Proof}}

\documentclass[12pt]{article}
\usepackage{amsfonts}
\usepackage[all]{xy}
\title{On Zariski's pairs of $m$th canonical\\ discriminant curves
}\author{Vik.S. Kulikov
\thanks{Partly supported by RFFI
 (No. 96-01-00614) and INTAS (No. 96-0713).}}
\date{        }

\begin{document}
\maketitle
\begin{abstract}
In this note we give examples of Zariski's pairs $B_{1,m},\, B_{2,m}$ ($m\in \mathbb N$ and $m\geq 5$) of plane cuspidal curves such that \newline ($i$) $B_{i,m}$ is the discriminant curve of a generic morphism $f_{i,m}:S_i\to \mathbb P^2$, $i=1,\, 2$, 
\newline ($ii$) $S_1$ and $S_2$ are homeomorphic surfaces of general type, 
\newline  ($iii$) $f_{i,m}$ is given by linear three-dimensional subsystem of the $m$th canonical class of $S_i$.
\end{abstract}   
\section*{Introduction}

By definition, two plane curves $B_1,\,B_2\subset \mathbb P^2$ over $\mathbb C$ are called a Zariski pair if they have the same degree  and homeomorphic tubular neighborhoods in $\mathbb P^2$, but the pairs $(\mathbb P^2,B_1)$ and $(\mathbb P^2,B_2)$ are not homeomorphic. A set of plane curves $(B_1, ...,\, B_k)$ is called {\it Zariski's $k$-tuple} if for each $i\neq j$ the pair $(B_i,\, B_j)$ is a Zariski one.

Let $B\subset \mathbb P^2$ be an irreducible plane curve with ordinary cusps and nodes, as its only singularities. Denote by $2d$ the degree of $B$, and let $g$ be the genus of its desingularization, $c= \# \{ \mbox{cusps of}\, B\}$, and $n= \# \{ \mbox{nodes of} \, B\}$. The curve $B$ is called {\it the discriminant curve} of a generic morphism if there exists a finite morphism $f:S\to \mathbb P^2$, $\deg f\geq 3$,  
satisfying the following conditions:

$(i)$ $S$ is a non-singular irreducible projective surface;

$(ii)$ $f$ is unramified over $\mathbb P^2 \setminus B$; 

$(iii)$  $f^{*}(B)=2R+C$, where $R$ is irreducible and non-singular, and $C$ is reduced; 

$(iv)$  $f_{\mid R}:R\to B$ coincides with the 
normalization of $B$. \newline 
Such a morphism $f$ is called {\it  generic}.  

Two generic morphisms $(S_1,f_1)$, $(S_2,f_2)$ with the same discriminant curve $B$ are said to be equivalent if there exists an isomorphism $\varphi : S_1 \to S_2$ such that 
$f_1=f_2\circ \varphi $. 
The following assertion is known as Chisini's Conjecture. 
\begin{gip}
Let $B$ be the discriminant curve of a generic morphism $f:S\to \mathbb P^2$ of degree $\deg f \geq 5$. Then $f$ is uniquely determined by $B$.
\end{gip}

In \cite{Kul1}, Chisini's conjecture was proved for discriminant curves of almost all generic morphisms. In particular, Chisini's conjecture holds for the discriminant curves of the generic morphisms $f_m:S\to \mathbb P^2$ of the surfaces $S$ of general type, where $f_m$ is given by a linear three-dimensional subsystem of the $m$th canonical class of $S$, $m\in \mathbb N$. The discriminant curve $B$ of a generic morphism $f_m:S\to \mathbb P^2$, given by a three-dimensional subsystem of the $m$th canonical class of $S$, will be called {\it the $m$th canonical discriminant curve}.

Let $(\mathbb P^2,B_1)$ and $(\mathbb P^2,B_2)$ be two diffeomorphic (resp. homeomorphic) pairs. One can show that if $B_1$ is the discriminant curve of a generic morphism $(S_1,f_1)$, then 
$B_2$ is also the discriminant curve of some generic morphism $(S_2,f_2)$. Moreover, if  $(S_1,f_1)$ is uniquely determined by $B_1$, then the same is true for $(S_2,f_2)$ and  $S_1$ and $S_2$ are diffeomorphic (resp. homeomorphic). It is a natural question to ask (\cite{Kul1}):

{\bf Question}. {\it Let $S_1,S_2\subset \mathbb P^r$ be two homeomorphic (Case C)(resp. diffeomorphic (Case D)) surfaces of general type embedded by the $m$th canonical class, and let $B_i$ be the $m$th canonical discriminant curve of a generic projection of $S_i$ onto $\mathbb P^2$. Are pairs $(\mathbb P^2,B_1)$ and $(\mathbb P^2,B_2)$ homeomorphic (resp. diffeomorphic)?}

In case C, the answer to this Question is negative, as the following  theorem shows.
\begin{thm} For each positive integer $k$ there exists a countable set of Zariski's $k$-tuples $(B_{1,m},...,\, B_{k,m})$, $m\in \mathbb N$, $m\geq 5$, of $m$th canonical discriminant curves $B_{i,m}$ such that the corresponding generic morphisms $f_{i,m}:S_i\to \mathbb P^2$, $i=1,...,\, k$, are the morphisms of pairwise homeomorphic surfaces of general type $S_i$.
\end{thm}

The proof of this Theorem is based on technique developed in \cite{Kul1}. It also uses the examples constructed of homeomorphic (but not diffeomorphic) surfaces of general type with ample canonical class constructed by Catanese in \cite{Cat1} and \cite{Cat2}.
 
The Question in the case D remains open.

This paper was written during my stay at the Max-Planck-Institut f\"{u}r Mathematik in Bonn. It is a pleasure to thank the Institut for its hospitality and financial support.

\section{Proof of Theorem 1.}
In \cite{Cat1} and \cite{Cat2}, Catanese investigated smooth simple bidouble coverings $\varphi : S\to Q=\mathbb P^1\times \mathbb P^1$ of type $(a,b),(m,n)$; thus, $\varphi $ is a finite $(\mathbb Z/2)^2$ Galois covering branched along two generic curves of respective bidegrees $(2a,2b)$, $(2m,2n)$. Let $a,b,m,n$ be integers satisfying the following conditions:
$$
\begin{array}{ll}
a   > 2n\, , n\geq 3,\,  &  m >2b \, ,  b \geq 3\, ; \\
 a  \equiv n\, (mod \, 2),\, & b\equiv m\, (mod \, 2) .
\end{array}
$$
Put 
\begin{equation} \label{uvwz}
u=(n+a-2),\, \, v=(m+b-2),\, \, w=(a-n),\, \, z=(m-b).
\end{equation}
By \cite{Cat1}, for a bidouble covering $\varphi : S\to Q=\mathbb P^1\times \mathbb P^1$ of type $(a,b),(m,n)$, $S$ is a simply connected surface of general type whose canonical class is
\begin{equation}  K_S=\varphi ^{*}(\mathcal{O}_Q(u,v))
\end{equation}
and
\begin{equation} \label{K^2}
\displaystyle K_S^2=8uv,\, \, \chi=\chi (\mathcal{O}_S)=\frac{3}{2}uv+(u+v)+2-\frac{1}{2}wz.
\end{equation}
For $u,v,w,z$ satisfying (\ref{uvwz}), the canonical class $K_S$ is 2-divisible in $H^2(S,\mathbb Z)$ and the index  
$$ r=r(S)= \max \{ s\in \mathbb N\, \, |\, \, (1/s)K_S\in H^2(S,\mathbb Z)\} $$ 
is equal to the greatest common divisor $(u,v)$. 

Since $K_S$ is 2-divisible, the intersection form on $H_2(S,\mathbb Z)$ is even. Therefore, by Freedman's theorem (\cite{Fre}) two such bidouble coverings $S_1$ and $S_2$ of $Q$ are homeomorphic if and only if $K_{S_1}^2=K_{S_2}^2$ and  $\chi (\mathcal{O}_{S_1})=\chi (\mathcal{O}_{S_2})$. In \cite{Cat1} and \cite{Cat2}, it has been proved that for each integer $k$ there exists at least one $k$-tuple $S_1,...,\, S_k$ of bidouble coverings of $Q$ of respective types $(a_i,b_i),(m_i,n_i)$ satisfying the conditions described above and such that 

$(i)$ $S_i$ and $S_j$ are homeomorphic for each $1\leq i,j\leq k $.

$(ii)$ $r(S_i)\neq r(S_j)$ for $i\neq j$ (and therefore $S_i$ and $S_j$ are not diffeomorphic). \newline We shall call such $k$-tuple of bidouble coverings of $Q$ {\it a Catanese $k$-tuple}. \vspace{0.3cm}
\newline {\bf Example.} Two bidouble coverings $S_1$ of type $(16,22),(52,4)$ and $S_2$ of type $(28,10),(28,10)$ form a Catanese pair with $$K_S^2=K_{S_i}^2=10368, \hspace{1cm} \chi=\chi (\mathcal{O}_{S_i})=1456$$ and 
$$r(S_1)=18, \hspace{1cm} r(S_2)=36.$$  \vspace{0.1cm} 

Consider the $m$-canonical embedding $\Phi _{|mK_{S_i}|}:S_i\to \mathbb P^r$, $m\geq 5$ (\cite{Bom}), and let $f_{m,i}:S_i \to \mathbb P^2$ be the restriction to $\Phi _{|mK_{S_i}|}(S_i)$ of a generic projection $pr:\mathbb P^r \to \mathbb P^2$. Denote by $R_{m,i}\subset S_i$ the ramification curve of the generic projection $f_{m,i}$, and by $B_{m,i}\subset \mathbb P^2$ the discriminant curve of $f_{m,i}$.

\begin{pred} \label{pred} Let $f_i:S_i\to \mathbb P^2$, $i=1,2$, be two generic morphisms with respective discriminant curves $B_i$. Assume that there exist a homeomorphism of pairs 
$\psi :(\mathbb P^2, B_1)\to (\mathbb P^2, B_2)$. Then there exists a 
commutative diagram 
$$
\xymatrix{
S_1 \ar[d]_{f_1} \ar[r]^{\Psi } & S_2 \ar[d]^{f_2} \\
\mathbb P^2 \ar[r]_{\psi } & \mathbb P^2 }
$$
such that 

($i$) $\Psi $ is also a homeomorphism;

($ii)$ $\Psi (R_1)=R_2$, where $R_i$ is the branch curve of $f_i$.
\end{pred}
\dvo  {} coincides with the proof of proposition 8 in \cite{Kul1}.

Let a generic morphism $f_m:S\to \mathbb P^2$ be given by a three-dimensional linear subsystem of the $m$th canonical class on $S$. Then, by the adjunction formula, the ramification curve $R_m\subset S$ of $f_m$  is equivalent to $(3m+1)K_S$. Therefore, if we have two generic morphisms $f_{i,m}:S_i\to \mathbb P^2$ given by three-dimensional linear subsystems of the $m$th canonical classes and such that the corresponding pairs $(\mathbb P^2,B_{1,m})$ and  $(\mathbb P^2,B_{2,m})$ are homeomorphic, then by proposition \ref{pred} there exists a homeomorphism $\Psi: S_1\to S_2$ such that $\Psi (R_{1,m})=R_{2,m}$. In this case, we should have $r(S_1)=r(S_2)$, since $\Psi $ induces an isomorphism $\Psi _{*}:H_2(S_1,\mathbb Z)\to H_2(S_2,\mathbb Z)$ such that $\Psi _{*}(R_{1,m})=R_{2,m}$. But it is impossible, since, by definition, $r(S_1)\neq r(S_2)$ for the Catanese pair $S_1,\, S_2$. Theorem 1 is proved. \vspace{0.3cm}
\newline {\bf Example.} (Continuation). Applying Theorem 1 for $m$-canonical discriminant curves of generic morphisms of  surfaces $S_1$ and $S_2$ in the example considered above, we obtain Zariski's pairs $(B_{1,m}, B_{2,m})$, $m\geq 5$, of plane cuspidal curves. Applying the computation in the proof of theorem 3 in \cite{Kul1}, one can see that the $m$-canonical discriminant curve $B_{i,m}$ has degree 
$$\deg B_{i,m}=10368m(3m+1),$$ genus $$g=5184(3m+2)(3m+1)+1,$$ and 
$$c=10368(12m^2+9m)-13632$$ ordinary cusps.

Steklov Mathematical Institute

{\rm victor$@$olya.ips.ras.ru }

\end{document}